\documentclass[11pt,reqno]{amsart}
\usepackage{epic}

\newtheorem{theorem}{Theorem}[]
\newtheorem{lemma}[theorem]{Lemma}

\newcommand\qbin[2]{{\left[\begin{matrix} #1 \\ #2 \end{matrix} \right]}}

\begin{document}

\title[Identities associated with mock theta functions $\omega(q)$ and $\nu(q)$ ]
{Some Identities associated with mock theta functions $\omega(q)$ and $\nu(q)$ }

\author{George E. Andrews}
\address{Department of Mathematics, The Pennsylvania State University,
University Park, PA 16802, USA} \email{gea1@psu.edu}

\author{Ae Ja Yee}
\address{Department of Mathematics, The Pennsylvania State University,
University Park, PA 16802, USA} \email{auy2@psu.edu}


\maketitle

\footnotetext[1]{The second author was partially supported by a grant ($\#$280903) from the Simons Foundation.} \vspace{0.5in}
\footnotetext[2]{2010 AMS Classification Numbers: Primary, 05A17; Secondary, 11P81.}

\noindent{\footnotesize{\bf Abstract.}  }  Recently, Andrews, Dixit and Yee defined two partition functions $p_{\omega}(n)$ and $p_{\nu}(n)$  that are related with Ramanujan's mock theta functions $\omega(q)$ and $\nu(q)$, respectively.  In this paper, we present two variable generalizations of their results. As an application, we reprove their results on $p_{\omega}(n)$ and $p_{\nu}(n)$ that are analogous to Euler's pentagonal number theorem. 

\medskip


\section{Introduction}



In \cite[p.~62]{watson0}, Watson defined the third order mock theta functions $\omega(q)$ and $\nu(q)$ as
\begin{align*}
\omega(q)&:=\sum_{n=0}^{\infty} \frac{q^{2n^2+2n}}{(q;q^2)_{n+1}^2}, 
\quad \nu(q):=\sum_{n=0}^{\infty} \frac{q^{n^2+n}}{(-q;q^2)_{n+1}},
\end{align*}
which were  rediscovered in Ramanujan's lost notebook \cite{geabcblnb5}, \cite[p.~15, p.~31]{lnb}.

Throughout the paper, we adopt the following $q$-series notation: 
\begin{align*}
(a;q)_0&:=1,\\
(a;q)_n&:=(1-a)(1-aq)\cdots (1-aq^{n-1}) \quad \text{for $n\ge 1$},\\
(a;q)_{\infty}&:=\lim_{n\to \infty} (a;q)_n, \quad |q|<1.
\end{align*}

The functions $\omega(q)$ and $\nu(q)$ along with other classical mock theta functions are well studied by basic hypergeometric series.  Andrews \cite{gea5} introduced the following two variable generalizations:
\begin{align}
\omega(z;q):=\sum_{n=0}^{\infty} \frac{z^n q^{2n^2+2n}}{(q;q^2)_{n+1}(zq;q^2)_{n+1}}, \quad 
\nu(z; q):=\sum_{n=0}^{\infty} \frac{ q^{n^2+n}}{(-zq;q^2)_{n+1}}. \label{andrews1}
\end{align}
Note that $\omega(1;q)=\omega(q)$ and $\nu(1,q)=\nu(q)$. In a subsequent paper \cite{gea6}, Andrews showed that
\begin{align}
\omega(z;q)=\sum_{n=0}^{\infty} \frac{z^n q^{n}}{(q;q^2)_{n+1}}, \quad 
\nu(z;q)=\sum_{n=0}^{\infty} (q/z;q^2)_n (-zq)^n. \label{andrews2}
\end{align}
These $\omega(z;q)$ and $\nu(z; q)$ are generalized further by Choi \cite{choi}.


Recently \cite{ady1}, Andrews Dixit, and Yee introduced two partition functions $p_{\omega}(n)$ and $p_{\nu}(n)$, where 
 $p_{\omega}(n)$ counts the number of partitions of $n$ in which  all odd parts are less than twice the smallest part, and  $p_{\nu}(n)$ counts the number of partitions of $n$ with the same constraints as $p_{\omega}(n)$ plus all parts being distinct.  Surprisingly, it was shown that these functions are very closely related with  $\omega(q)$ and $\nu(q)$. Namely, 
\begin{align}
\sum_{n=1}^{\infty} p_{\omega}(n) &=q\, \omega(q),\quad 
\sum_{n=0}^{\infty} p_{\nu}(n)= \nu(-q), \label{omeganu}
\end{align}
By the definitions of $p_{\omega}(n)$ and $p_{\nu}(n)$, \eqref{omeganu} can be stated as follows:
\begin{align}
\sum_{n=1}^{\infty} \frac{q^n}{(q^{n};q)_{n+1} (q^{2n+2};q^2)_{\infty}}&=\sum_{n=0}^{\infty} \frac{q^{2n^2+2n+1}}{(q;q^2)_{n+1}}, \label{omega}\\
\sum_{n=0}^{\infty} q^n (-q^{n+1};q)_n (-q^{2n+2};q^2)_{\infty}&=\sum_{n=0}^{\infty} \frac{q^{n^2+n}}{(q;q^2)_{n+1}}. \label{nu}
\end{align}

The main purpose of this paper is to provide two variable generalizations of  \eqref{omega} and \eqref{nu}, which are given in the following theorem. 

\begin{theorem}\label{thm1}
We have
\begin{align}
\sum_{n=1}^{\infty} \frac{q^n}{(zq^n;q)_{n+1} (zq^{2n+2};q^2)_{\infty}}  
&=\sum_{n=0}^{\infty} \frac{z^nq^{2n^2+2n+1}}{(q;q^2)_{n+1} (zq;q^2)_{n+1}}
,\label{eqthm1}
\\
\sum_{n=0}^{\infty} q^n(-zq^{n+1};q)_{n} (-zq^{2n+2};q^2)_{\infty} &=\sum_{n=0}^{\infty} \frac{z^n q^{n^2+n}}{(q;q^2)_{n+1}}.\label{eqthm2}
\end{align}
\end{theorem}

It is immediate that \eqref{eqthm1} and \eqref{eqthm2} yield \eqref{omega} and \eqref{nu}, respectively, when $z=1$.  We also note that by \eqref{andrews1}, \eqref{andrews2}, and \eqref{eqthm1}, we have
\begin{equation}
\sum_{n=1}^{\infty} \frac{q^{n}}{(zq^n;q)_{n+1} (zq^{2n+2};q^2)_{\infty}}=\sum_{n=1}^{\infty}  \frac{z^{n-1} q^{n}}{(q;q^2)_n}. \label{thm2omega}
\end{equation}

Meanwhile, the sum on the right hand side of \eqref{eqthm2} with $q$ replaced by $-q$ does not match $\nu(z;q)$ in \eqref{andrews1}.  Thus, Theorem~\ref{thm1} suggests a different two variable generalization of $\nu(q)$. 
Let
\begin{equation*}
\nu_1(z;q):=\sum_{n=0}^{\infty} \frac{z^n q^{n^2+n}}{(-q;q^2)_{n+1}}.
\end{equation*}
Then, $\nu_1(z;q)$ also has a representation that is reminiscent to that of $\nu(z;q)$ in  \eqref{andrews2}. 
\begin{theorem} \label{thma}
We have
\begin{equation}
\nu_1(z;q)=\sum_{n=0}^{\infty} (zq;q^2)_n (-q)^n. \label{B}
\end{equation}
\end{theorem}

In \cite{ady1}, analogues of Euler's pentagonal number theorem for $p_{\omega}(n)$ and $p_{\nu}(n)$ were given as follows:
\begin{align}
\sum_{n=1}^{\infty} \frac{q^n}{(-q^n;q)_{n+1}(-q^{2n+2};q^2)_{\infty}} &=\sum_{j=0}^{\infty} (-1)^j q^{6j^2+4j+1}(1+q^{4j+2}), \label{omegaeuler}\\
\sum_{n=0}^{\infty} q^n (q^{n+1};q)_n (q^{2n+2};q^2)_{\infty} &=\sum_{j=0 }^{\infty} (-1)^j q^{3j^2+2j} (1+q^{2j+1}). \label{nueuler}
\end{align}
These analogues will be reproved using Theorem~\ref{thm1}  with some identities from Ramanujan's lost notebook \cite{geabcbln1}.
 

The proof of \eqref{omeganu} given in \cite{ady1} uses a complex four-parameter $q$-series identity \cite[p. 141, Theorem 1]{gea90}, whereas surprisingly, the bulk of the proof of Theorem~\ref{thm1} relies on two finite summations (see Lemmas~\ref{lem4} and \ref{lem5}), namely
 \begin{equation*}
 \sum_{j=0}^{n} \frac{q^j (q;q)_{n+j}}{(q^2;q)_j} =(q^2;q^2)_n
 \end{equation*}
 and
 \begin{equation*}
 \sum_{j=0}^{n} \frac{q^{2j}(q;q)_{n+j}}{(q^2;q^2)_j}=(q;q^2)_{n+1}+q^{n+1}(q^2;q^2)_n.
 \end{equation*}
 These results, while seemingly quite simple, are surprising for a couple of reasons. First, the sums do not terminate naturally, and second, we were unable to find these results in the $q$-series literature.

The rest of this paper is organized as follows. In Section~\ref{seca}, Theorem~\ref{thma} will be proved.  In Section~\ref{sec2}, we provide necessary lemmas for the proof of Theorem~\ref{thm1}, 
and then we will prove Theorem~\ref{thm1} in Section~\ref{sec3}. Finally, \eqref{omegaeuler} and \eqref{nueuler} will be proved in Section~\ref{sec4}.  We conclude our paper with some remarks in Section~\ref{sec5}.

\section{Two variable generalizations $\omega(z,q)$ and $\nu_1(z,q)$} \label{seca}

We first prove a theorem on $\omega(z,q)$.
\begin{theorem}  \label{lem1}
We have
\begin{align}
\omega(z;q)&
=\sum_{n=0}^{\infty} \frac{q^n}{(zq;q^2)_{n+1}}.\label{eqlem1} 
\end{align}
\end{theorem}

\proof 
By \eqref{andrews2}, 
it will be sufficient to show that
\begin{equation*}
\sum_{n=0}^{\infty} \frac{z^{n} q^n}{(q;q^2)_{n+1}}=\sum_{n=0}^{\infty} \frac{q^n}{(zq;q^2)_{n+1}}..
\end{equation*}
Recall the $q$-binomial coefficient $\qbin{n}{m}_q$ defined by \cite[p.~35]{gea1998}
\begin{equation*}
\qbin{n}{m}_{q}:=
\begin{cases}
\frac{(q;q)_{n}}{(q;q)_m(q;q)_{n-m}},  & 0\leq m\leq n,\\
0, & \text{otherwise}.
\end{cases}
\end{equation*}

Using the $q$-binomial theorem  \cite[p. 36, Equation (3.3.7)]{gea1998},
\begin{align*}
\sum_{n=0}^{\infty} \frac{z^nq^n}{(q;q^2)_{n+1}}&=
\sum_{n=0}^{\infty} z^nq^{n} \sum_{m=0}^{\infty} q^{m} \qbin{n+m}{m}_{q^2}\\
&=\sum_{m= 0}^{\infty}  q^{m} \sum_{n=0}^{\infty} z^n q^{n} \qbin{n+m}{m}_{q^2}\\
&=\sum_{m= 0}^{\infty} \frac{ q^{m}}{(zq;q^2)_{m+1}}.
\end{align*}
\endproof

We now prove Theorem~\ref{thma}.

\noindent{\it Proof of Theorem~\ref{thma}}. 
Again, we use the $q$-binomial theorem  \cite[p. 36, Equation (3.3.7)]{gea1998} to see
\begin{align*}
\sum_{n=0}^{\infty} \frac{z^n q^{n^2+n} }{(-q;q^2)_{n+1}}&=\sum_{n=0}^{\infty} z^n q^{n^2+n} \sum_{m=0}^{\infty} (-q)^m \qbin{n+m}{m}_{q^2}\\
&= \sum_{N=0}^{\infty} (-q)^N \sum_{n=0}^{N} (-1)^n z^n q^{n^2} \qbin{N}{n}_{q^2}\\
&=\sum_{N=0}^{\infty}  (-q)^N ( zq;q^2)_N,
\end{align*}
where for the last equality, the $q$-binomial theorem  \cite[p. 36, Equation (3.3.6)]{gea1998} is used. 
{\hfill $\Box$ \par\medskip}

\medskip

\noindent{\bf Remarks.} 1. In the Rogers-Fine identity \cite[Chap.9, p. 223]{geabcbln1}, send $q\to q^2$, then $\alpha=0$, $\beta=q^3$, and multiply both sides by $q/(1-q)$. The result is the first identity on $\omega(z;q)$ in \eqref{andrews2}.

2. In \cite[p. 29, Ex. 6]{gea1998}, set $x=zq$ and $y=-q$. Then the result is  Theorem~\ref{thma}. Alternatively,  
in the second Heine's transformation:
\begin{equation*}
_{2} \phi_1 \left(\begin{matrix} a, & b\\ & c\end{matrix} \, ; q, t\right)=\frac{(c/b)_{\infty} (bt)_{\infty}}{(c)_{\infty} (t)_{\infty}} {_2} \phi_1 \left(\begin{matrix} abt/c, & b\\ & bt \end{matrix}\, ; q, c/b\right),
\end{equation*}
replace $q$ by $q^2$, set $c=-q^3, b=q^2, a=-q^2z/t$ and let $t\to 0$. Then multiply both sides by $1/(1+q)$. Here,
\begin{equation*}
_{2} \phi_1 \left(\begin{matrix} a, & b\\ & c\end{matrix} \, ; q, t\right):=\sum_{n=0}^{\infty} \frac{(a;q)_n (b;q)_n}{(q;q)_n (c;q)_n} t^n. 
\end{equation*}

 
 
 
 
\section{Lemmas} \label{sec2}

For any $n, i\ge 0$,  let 
\begin{equation*}
S_n(i):=\sum_{j=0}^{n} \frac{q^{ij} (q;q)_{n+j}}{(q^2;q^2)_j}.
\end{equation*}

We need some functional equations for $S_n(i)$. 
\begin{lemma} \label{lem2}
We have
\begin{equation*}
S_n(i)=S_{n-1}(i) -q^n S_{n-1}(i+1) +q^{in} (q;q^2)_n.
\end{equation*}
\end{lemma}

\proof
\begin{align*}
S_n(i)-S_{n-1}(i)&=\sum_{j=0}^{n-1} \frac{q^{ij} (q;q)_{n-1+j} (1-q^{n+j} -1)}{(q^2;q^2)_j} +\frac{q^{in}(q;q)_{2n}}{(q^2;q^2)_n}\\
&=-q^n S_{n-1}(i+1)+q^{in} (q;q^2)_n.
\end{align*}
\endproof

\begin{lemma} \label{lem3}
We have
\begin{equation*}
S_n(i+2)=S_{n}(i)-q^{i} S_{n+1}(i)+ q^{i(n+1)}(1+q^i)(q;q^2)_{n+1}.
\end{equation*}
\end{lemma}

\proof
\begin{align*}
S_n(i+2) &=\sum_{j=0}^{n} \frac{q^{ij} \big(1-(1-q^{2j})\big)(q;q)_{n+j}}{(q^2;q^2)_j} \\
&=S_{n}(i)-\sum_{j=1}^n \frac{q^{ij} (q;q)_{n+j}}{(q^2;q^2)_{j-1}}\\
&=S_n(i)-\sum_{j=0}^{n-1} \frac{q^{i(j+1)} (q;q)_{n+1+j}}{(q^2;q^2)_{j}}\\
&=S_n(i)-q^i S_{n+1}(i) +\frac{q^{i(n+1)} (q;q)_{2n+1}}{(q^2;q^2)_n}+ \frac{q^{i(n+2)}(q;q)_{2n+2}}{(q^2;q^2)_{n+1}}\\
&=S_n(i)-q^i S_{n+1}(i)+ q^{i(n+1)}(1+q^i)(q;q^2)_{n+1}.
\end{align*}
\endproof

Using the previous lemmas, we evaluate $S_n(i)$ for $i=1,2$. 
\begin{lemma}  \label{lem4}
We have
\begin{equation*}
S_n(1)=(q^2;q^2)_n.
\end{equation*}
\end{lemma}

\proof
By Lemma~\ref{lem2}, 
\begin{align}
S_n(1)&=S_{n-1}(1) -q^n S_{n-1}(2) +q^n (q;q^2)_n, \label{eq1}\\
S_{n-1}(2)&=S_{n-2}(2)  -q^{n-1} S_{n-2}(3) +q^{2n-2} (q;q^2)_{n-1}, \label{eq2}
\end{align}
and by Lemma~\ref{lem3} with $i=1$ and $n$ replaced by $n-2$, 
\begin{align}
S_{n-2}(3)&=S_{n-2}(1)-qS_{n-1}(1)+ q^{n-1}(1+q) (q;q^2)_{n-1}. \label{eq3}
\end{align}
Now by \eqref{eq1}, 
\begin{align}
S_{n-1}(2)&=q^{-n} \big(-S_{n}(1)+S_{n-1}(1)+ q^{n}(q;q^2)_{n} \big). \label{eq4}
\end{align}

We now obtain a recurrence for $S_n(1)$ by substituting from \eqref{eq3} and \eqref{eq4} into \eqref{eq2}. Hence
\begin{align}
&q^{-n} \big(-S_n(1)+S_{n-1}(1)+ q^n (q;q^2)_n\big) \notag\\
&=q^{-n+1}\big(-S_{n-1}(1)+S_{n-2}(1)+q^{n-1} (q;q^2)_{n-1}\big)\notag\\
& -q^{n-1}\big(S_{n-2}(1)-q S_{n-1}(1)+q^{n-1} (1+q)(q;q^2)_{n-1}\big)  +q^{2n-2} (q;q^2)_{n-1}. \label{eq5}
\end{align}
We now simplify \eqref{eq5} to obtain 
\begin{align*}
S_n(1)=(1+q-q^{2n})S_{n-1}(1)-q(1-q^{2n-2})S_{n-2}(1). 
\end{align*}
However, we see that
\begin{align*}
&(1+q-q^{2n})(q^2;q^2)_{n-1} - q(1-q^{2n-2})(q^2;q^2)_{n-2}\notag\\
&=(q^2;q^2)_n +q(q^2;q^2)_{n-1} -q(q^2;q^2)_{n-1}\\
&=(q^2;q^2)_n.
\end{align*}
Thus both $S_n(1)$ and $(q^2;q^2)_n$ satisfy the same recurrence and
\begin{align*}
S_0(1)&=1=(q^2;q^2)_0,\\
S_1(1)&=1-q^2=(q^2;q^2)_1.
\end{align*}
Therefore, by mathematical induction,
\begin{align*}
S_n(1)=(q^2;q^2)_n. 
\end{align*}
for every $n$. 
\endproof

\begin{lemma} \label{lem5}
We have
\begin{equation*}
S_n(2)=(q;q^2)_{n+1}+q^{n+1}(q^2;q^2)_n.
\end{equation*}
\end{lemma}

\proof
By Lemma~\ref{lem2}, with $i=1$, $n$ replaced by $n+1$,
\begin{align*}
S_n(2)&=q^{-n-1}\big(S_n(1)-S_{n+1}(1)+q^{n+1} (q;q^2)_{n+1}\big)\\
&=q^{-n-1} \big((q^2;q^2)_n -(q^2;q^2)_{n+1} \big) +(q;q^2)_{n+1} \tag{by Lemma~\ref{lem4}}
\\
&=q^{-n-1}(q^2;q^2)_n (1-1+q^{2n+2}) +(q;q^2)_{n+1}\\
&=q^{n+1} (q^2;q^2)_n+(q;q^2)_{n+1}.
\end{align*}
\endproof

\section{Proof of Theorem~\ref{thm1}}\label{sec3}

\subsection{Proof of \eqref{eqthm1}}
We will prove the equivalent identity \eqref{thm2omega}.  First, we expand the two products in the denominator on the left hand side of \eqref{thm2omega} into power series in $z$ and then we compare coefficients of $z^N$.  Hence \eqref{thm2omega} is equivalent to proving
\begin{equation*}
\sum_{n=1}^{\infty} \sum_{s=0}^{N} q^{n(1+N+s)+2s}\qbin{n+N-s}{n} \frac{1}{(q^2;q^2)_s}=\frac{q^{N+1}}{(q;q^2)_{N+1}},
\end{equation*}
and summing the $n$ series by the $q$-binomial theorem \cite[p. 36, Equation (3.3.7)]{gea1998}, we find that \eqref{thm2omega} is equivalent to the assertion that 
\begin{equation}
\sum_{s=0}^N \frac{q^{2s}}{(q^2;q^2)_s} \bigg(\frac{1}{(q^{1+N+s};q)_{N-s+1}}-1\bigg)=\frac{q^{N+1}}{(q;q^2)_{N+1}}. \label{eq111}
\end{equation}

Now mathematical induction reveals immediately that 
\begin{equation*}
\sum_{s=0}^{N} \frac{q^{2s}}{(q^2;q^2)_s} =\frac{1}{(q^2;q^2)_N}.
\end{equation*}
Hence rewriting \eqref{eq111} using this fact, we see that \eqref{thm2omega} is equivalent to 
\begin{equation}
\sum_{s=0}^N \frac{q^{2s}}{(q^2;q^2)_s (q^{1+N+s};q)_{N-s+1}} =\frac{1}{(q^2;q^2)_N}+\frac{q^{N+1}}{(q;q^2)_{N+1}}. \label{eq12}
\end{equation}
Now multiply both sides by $(q;q)_{2N+1}$, and we obtain
\begin{equation}
\sum_{s=0}^N \frac{q^{2s} (q;q)_{N+s}}{(q^2;q^2)_s} =(q;q^2)_{N+1} + q^{N+1} (q^2;q^2)_N. \label{eq13}
\end{equation}
Since \eqref{eq13} is merely Lemma~\ref{lem5}, we see that \eqref{thm2omega}, consequently \eqref{eqthm1}, is proved. 

\subsection{Proof of \eqref{eqthm2}}

As in the proof of \eqref{eqthm1}, we expand the two products on the left hand side  of \eqref{eqthm2} and then compare the coefficients of $z^N$. Thus \eqref{eqthm2} is equivalent to
\begin{equation*}
\sum_{n= 0}^{\infty} \sum_{s=0}^{N} q^{n+\binom{N-s+1}{2}+n(N-s)+s^2+s+2ns} \qbin{n}{N-s} \frac{1}{(q^2;q^2)_s}=\frac{q^{N^2+N}}{(q;q^2)_{N+1}}, 
\end{equation*}
which is equivalent to 
\begin{equation}
\sum_{n=0}^{\infty} \sum_{s=0}^{N} q^{\binom{N-s+1}{2}+(n-N+s)(N+s+1)} \qbin{n}{n-N+s} \frac{1}{(q^2;q^2)_s}=\frac{1}{(q;q^2)_{N+1}}. \label{eq15}
\end{equation}

We now sum the series on $n$ by the $q$-binomial theorem \cite[p. 36, Equation (3.3.7)]{gea1998}. Hence \eqref{eq15} is equivalent to 
\begin{equation}
\sum_{s=0}^{N} \frac{q^{\binom{N-s+1}{2}}}{(q^2;q^2)_s (q^{N+s+1};q)_{N-s+1}}=\frac{1}{(q;q^2)_{N+1}}. \label{eq16}
\end{equation}

We now multiply \eqref{eq16} by $(q;q)_{2N+1}$ to obtain
\begin{equation}
\sum_{s=0}^{N} \frac{q^{\binom{N-s+1}{2}} (q;q)_{N+s}}{(q^2;q^2)_s} =(q^2;q^2)_N. \label{eq17}
\end{equation}

Now in \eqref{eq17} replace $q$ by $1/q$ and multiply by $q^{N^2+N}$. The resulting equivalent identity is
\begin{equation}
\sum_{s=0}^N \frac{q^s (q;q)_{N+s}}{(q^2;q^2)_s}=(q^2;q^2)_N. \label{eq18}
\end{equation}
But \eqref{eq18} is merely a restatement of Lemma~\ref{lem4}. Hence \eqref{eqthm2} is proved.

\section{Analogues of Euler's pentagonal number theorem for $p_{\omega}(n)$ and $p_{\nu}(n)$} \label{sec4}

In this section, we will reprove \eqref{omegaeuler} and \eqref{nueuler}.  We first prove \eqref{omegaeuler}. Setting $z=-1$ in \eqref{eqthm1} and \eqref{eqlem1}, we obtain
\begin{align}
\sum_{n=1}^{\infty} \frac{q^n}{(-q^n;q)_{n+1}(-q^{2n+2};q^2)_{\infty}}
&=\sum_{n=0}^{\infty} \frac{q^{n+1}}{(-q;q^2)_{n+1}}. \label{eq11}
\end{align}

Recall Entry 9.5.3 in Chapter 9, Lost Notebook Part 1 \cite{geabcbln1}: 
\begin{equation}
\sum_{n=0}^{\infty} \frac{q^n}{(-q;q^2)_{n+1}} = \sum_{j=0}^{\infty} (-1)^j q^{6j^2+4j} (1+q^{4j+2}). \label{eq12}
\end{equation}
Thus, \eqref{omegaeuler} follows from \eqref{eq11} and \eqref{eq12}. 

Similarly, for \eqref{nueuler},  set $z=-1$ in \eqref{eqthm2} and \eqref{B} to get
\begin{equation*}
\sum_{n=0}^{\infty} q^n (q^{n+1};q)_n (q^{2n+2};q^2)_{\infty}=\sum_{n=0}^{\infty} (q;q^2)_n q^n.
\end{equation*}
We then use Entry 9.5.2 in Chapter 9, Lost Notebook Part 1 \cite{geabcbln1}:
\begin{equation*}
\sum_{n=0}^{\infty} (q;q^2)_{n} q^n = \sum_{j=0}^{\infty} (-1)^j q^{3^2+2j} (1+q^{2j+1}).
\end{equation*}

\section{Concluding Remarks}\label{sec5}

The results in Theorem~\ref{thm1} cry out for a bijective proof. In each case, the insertion of $z$ into the sums refines the partitions being generated according to the number of parts in each partition. While it is not difficult at all to bijectively prove the identity of the simple series in each of \eqref{andrews2} and \eqref{B}, we have been completely unable to obtain bijectively the portion where the number of odds is restricted by the smallest parts.


\begin{thebibliography}{99}

\bibitem{gea5}
G. E.~Andrews, \emph{On basic hypergeometric series, mock theta functions, and partitions, I}, Quart. J. Math.~\textbf{17} (1966), 64--80. 

\bibitem{gea6}
G. E.~Andrews, \emph{On basic hypergeometric series, mock theta functions, and partitions, II}, Quart. J. Math.~\textbf{17} (1966), 132--143. 

\bibitem{gea1998}
G. E.~Andrews, The theory of partitions, Addison-Wesley Pub. Co., NY, 300 pp. (1976). Reissued, Cambridge University Press, New York, 1998.

\bibitem{gea90}
G. E.~Andrews, \emph{Ramanujan''s ``Lost" notebook: \textup{I}. Partial theta functions}, Adv. Math.~\textbf{41} (1981), 137--172.

\bibitem{geabcbln1}
G. E. Andrews and B.C. Berndt, Ramanujan's Lost Notebook, Part I. Springer, New York, 437 pp. (2005).

\bibitem{geabcblnb5}
G. E. Andrews and B.C. Berndt, Ramanujan's Lost Notebook, Part V., in preparation.

\bibitem{ady1}
G. E.~Andrews, A.~Dixit and A. J.~Yee, \emph{Partitions associated with the Ramanujan/Watson mock theta functions $\omega(q), \nu(q)$ and $\phi(q)$}, \emph{Research in Number Theory},~\textbf{1} Issue 1 (2015), 1--25.




\bibitem{choi}
Y.-S. Choi, \emph{The basic bilateral hypergeometric series and the mock theta functions}, Ramanujan J. \textbf{24} (2011), 345--386.


\bibitem{lnb}
S.~Ramanujan, \emph{The Lost Notebook and Other Unpublished Papers}, Narosa, New Delhi, 1988.


\bibitem{watson0}
G.N.~Watson, \emph{The final problem: an account of the mock theta functions}, J. London Math. Soc.~\textbf{11} (1936), 55--80.

\end{thebibliography}
\end{document}